\documentclass[12pt]{article}
\usepackage{amsmath,amssymb,amsthm}
\usepackage{hyperref}
\usepackage{graphicx}
\usepackage[mathlines]{lineno}
\usepackage{mathscinet}
\usepackage{dsfont} 
\usepackage{color}
\definecolor{red}{rgb}{1,0,0}

 \newtheorem{thm}{Theorem}[section]
 \newtheorem{defn}[thm]{Definition}
 
 \newtheorem{cor}[thm]{Corollary}
 \newtheorem{lem}[thm]{Lemma}
 
 \newtheorem{conj}[thm]{Conjecture}

 \newtheorem{obs}[thm]{Observation}
 
 \newtheorem{rem}[thm]{Remark}

\def\Krn#1#2{\bigvee^{#1} \overline{K_{ #2}}}

\DeclareMathOperator{\crn}{cr}
\newcommand{\rcrn}{\overline{\operatorname{cr}}}
\DeclareMathOperator{\CRN}{CR}

\newcommand{\oururl}{\url{https://orion.math.iastate.edu/lidicky/pub/knnn}}

\newcommand{\ff}{{\bf f}}

\newcommand{\bit}{\begin{itemize}}
\newcommand{\eit}{\end{itemize}}
\newcommand{\ben}{\begin{enumerate}}
\newcommand{\een}{\end{enumerate}}
\newcommand{\beq}{\begin{equation}}
\newcommand{\eeq}{\end{equation}}
\newcommand{\bea}{\begin{eqnarray*}}
\newcommand{\eea}{\end{eqnarray*}}
\newcommand{\bean}{\begin{eqnarray}}
\newcommand{\eean}{\end{eqnarray}}
\newcommand{\bpf}{\begin{proof}}
\newcommand{\epf}{\end{proof}\ms}
\newcommand{\bmt}{\begin{bmatrix}}
\newcommand{\emt}{\end{bmatrix}}
\newcommand{\ms}{\medskip}
\newcommand{\noi}{\noindent}
\newcommand{\ba}{\begin{array}}
\newcommand{\ea}{\end{array}}

\begin{document}
\title{Crossing numbers of complete tripartite and balanced complete multipartite graphs}
\author{Ellen Gethner\thanks{Department of Computer Science and Engineering, University of Colorado Denver,
(ellen.gethner@ucdenver.edu).}, Leslie Hogben\thanks{Department of Mathematics, Iowa State University 
and American Institute of
Mathematics 
(hogben@aimath.org).}, Bernard Lidick\'y\thanks{Department of Mathematics, Iowa State University 
(lidicky@iastate.edu, myoung@iastate.edu).}, 
Florian Pfender\thanks{Department of Mathematical and Statistical Sciences, University of Colorado Denver,
(Florian.Pfender@ucdenver.edu).}, \\Amanda Ruiz\thanks{Department of Mathematics and Computer Science, University of San Diego (alruiz@sandiego.edu).}, Michael Young\footnotemark[3]}


\maketitle

\begin{abstract}
The crossing number $\crn(G)$ of a graph $G$ is the minimum number
of crossings in a nondegenerate planar drawing of $G$. 
The rectilinear crossing
number $\overline{\operatorname{cr}}(G)$  of $G$ is the minimum number of crossings in a rectilinear nondegenerate planar drawing (with edges as straight line segments) of
$G$.   
Zarankiewicz proved in 1952 that  
 $\overline{\operatorname{cr}}(K_{n_1,n_2})\le Z(n_1,n_2):=\left\lfloor \frac {n_1}{2}\right\rfloor\left\lfloor \frac {n_1-1} 2\right\rfloor\left\lfloor \frac {n_2} 2\right\rfloor\left\lfloor \frac {n_2-1} 2\right\rfloor$.  We define an analogous bound for the complete tripartite graph $K_{n_1,n_2,n_3}$,
 \[A(n_1,n_2,n_3)  =\!\!\!\!\!\sum_{i=1,2,3\atop \{j,k\}=\{1,2,3\}\setminus \{i\}}\!\!\!\!\!\left(
\left\lfloor \frac{n_j\vphantom{1}}{2}\right\rfloor \left\lfloor \frac{n_j-1}{2}\right\rfloor  
\left\lfloor \frac{n_k\vphantom{1}}{2}\right\rfloor \left\lfloor \frac{n_k-1}{2}\right\rfloor +
\left\lfloor \frac{n_i\vphantom{1}}{2}\right\rfloor \left\lfloor \frac{n_i-1}{2}\right\rfloor \left\lfloor \frac{n_j n_k\vphantom{1}}{2}\right\rfloor \right)\!\!,\] and prove $\rcrn(K_{n_1,n_2,n_3})\le A({n_1,n_2,n_3})$.  We also show that for $n$ large enough,  $0.973 A(n,n,n) \le \rcrn(K_{n,n,n})$ and $0.666 A(n,n,n)\le \crn(K_{n,n,n})$, with the tighter rectilinear lower bound established through the use of flag algebras.

A  complete multipartite graph is  balanced if the partite sets all have the same cardinality.
We study asymptotic behavior of the crossing number of the balanced complete $r$-partite graph. Richter and Thomassen proved in 1997 that the limit as $n\to\infty$ of $\crn(K_{n,n})$ over the maximum number of crossings in a drawing of $K_{n,n}$ exists and is at most $\frac 1 4$.  We define
$\zeta(r)=\frac{3(r^2-r) }{8 \left(r^2+r-3\right)}$ and show that for a fixed $r$ and the balanced complete $r$-partite graph, $\zeta(r)$ is an upper bound to the limit superior of the crossing number divided by the maximum number of crossings in a drawing.

\end{abstract}

 \medskip



\section{Introduction}\label{sintro}

This paper deals with two main topics, the crossing numbers and rectilinear crossing
numbers of complete tripartite graphs, and the asymptotic behavior of the crossing number of a balanced complete multipartite graph. 
   In  the introduction,  we provide background, present definitions,  state our main results, and make related conjectures.  

A plane drawing of a graph is a 
{\em good drawing} if no more than two edges intersect at any point that is not a vertex,  edges incident with a common vertex do not cross,  no pair of
edges cross more than once, and edges that intersect  at a non-vertex must cross.  The {\em crossing number} of a good drawing $D$ is the number of non-vertex edge intersections in $D$.  The {\em crossing number} of a graph $G$ is 
\[\crn(G):=\min\{\crn(D) : D \mbox{ is a good drawing of }G\}.\]
Clearly a graph $G$ is planar if and only if $\crn(G)=0$. Tur\' an contemplated the question of determining the crossing number of the complete bipartite
graph $K_{n,m}$ during World War II, as described in \cite{Turan}.   After he posed the problem in lectures in Poland in 1952,  Zarankiewicz \cite{Zarank} proved that  
  \[\crn(K_{n,m})\le Z(n,m):=\left\lfloor \frac {n\vphantom{1}} 2\right\rfloor\left\lfloor \frac {n-1} 2\right\rfloor \left\lfloor \frac {m\vphantom{1}} 2\right\rfloor\left\lfloor \frac {m-1} 2\right\rfloor\]   
 and attempted to prove $\crn(K_{n,m})= Z(n,m)$; the latter equality has become known as Zarankiewicz's Conjecture. Hill's Conjecture for the crossing number of the complete graph $K_n$ is 
  \[\crn(K_{n})=H(n):=\frac 1 4 \left\lfloor \frac {n\vphantom{1}} 2\right\rfloor \left\lfloor \frac {n-1} 2\right\rfloor \left\lfloor \frac {n-2} 2\right\rfloor \left\lfloor \frac {n-3} 2\right\rfloor\!,\] 
and it is known  that $\crn(K_{n})\le H(n)$.   Background on crossing numbers, including these well-known conjectures, can be found in  \cite{BW10} and \cite{RT97MAA}.

We establish an  upper bound for the rectilinear crossing number of a complete tripartite graph that is analogous to Zarankiewicz' bound.   Define 
\[A(n_1,n_2,n_3) :=\!\!\!\!\!\sum_{i=1,2,3\atop \{j,k\}=\{1,2,3\}\setminus \{i\}}\!\!\!\left(
\left\lfloor \frac{n_j\vphantom{1}}{2}\right\rfloor \left\lfloor \frac{n_j-1}{2}\right\rfloor  
\left\lfloor \frac{n_k\vphantom{1}}{2}\right\rfloor \left\lfloor \frac{n_k-1}{2}\right\rfloor +
\left\lfloor \frac{n_i\vphantom{1}}{2}\right\rfloor \left\lfloor \frac{n_i-1}{2}\right\rfloor \left\lfloor \frac{n_j n_k\vphantom{1}}{2}\right\rfloor \right)\!.\] 

 Very little is known about exact values of crossing numbers of complete tripartite graphs, except when two of the parts are small. 
  For example, $\crn(K_{1,3,n})=2\left\lfloor \frac{n}{2}\right\rfloor \left\lfloor \frac{n-1}{2}\right\rfloor + \left\lfloor \frac{n}{2}\right\rfloor$ and $\crn(K_{2,3,n})=4\left\lfloor \frac{n}{2}\right\rfloor \left\lfloor \frac{n-1}{2}\right\rfloor + n$ are established in \cite{A86}, $\crn(K_{1,4,n})=n(n-1)$  is established in \cite{Ho08, HZ08},  and 
$\crn(K_{2,4,n})=6\left\lfloor \frac{n}{2}\right\rfloor \left\lfloor \frac{n-1}{2}\right\rfloor + 2n$  is established in \cite{Ho13}.  
  It is straightforward to verify that $A(1,3,n)=2\left\lfloor \frac{n}{2}\right\rfloor \left\lfloor \frac{n-1}{2}\right\rfloor + \left\lfloor \frac{n}{2}\right\rfloor=\crn(K_{1,3,n})$, $A(2,3,n)=4\left\lfloor \frac{n}{2}\right\rfloor \left\lfloor \frac{n-1}{2}\right\rfloor + n=\crn(K_{2,3,n})$,  $A(1,4,n)=n(n-1)=\crn(K_{1,4,n})$  and 
$A(2,4,n)=6\left\lfloor \frac{n}{2}\right\rfloor \left\lfloor \frac{n-1}{2}\right\rfloor + 2n=\crn(K_{2,4,n})$.

A good planar drawing of $G$ is {\em rectilinear} if every edge is
drawn as a straight line segment, and the {\em rectilinear crossing
number} $\rcrn(G)$  of $G$ is the minimum number of crossings in a rectilinear  drawing of
$G$; clearly $\crn(G)\le \rcrn(G)$. Zarankiewicz proved  that 
$\crn(K_{n,m})\le Z(n,m)$ by exhibiting a drawing that actually proves $\rcrn(K_{n,m})\le Z(n,m)$, because the drawing  is rectilinear.

The next three theorems give bounds on the crossing number and rectilinear crossing number of complete tripartite graphs and are proved in  Section \ref{sA3n}.

 \begin{thm} \label{corAn1n2n3}
For all $n_1,n_2,n_3\ge 1$,   $\crn(K_{n_1,n_2,n_3})\le\rcrn(K_{n_1,n_2,n_3})\le A(n_1,n_2,n_3)$.
\end{thm}

 \begin{thm}\label{LBthm} For $n$ large enough,  $0.666 A(n,n,n) \le \crn(K_{n,n,n})$.
\end{thm}


 \begin{thm}\label{rectLBthm} For $n$ large enough,  $0.973 A(n,n,n) \le \rcrn(K_{n,n,n})$.
\end{thm}

Theorem \ref{LBthm} is proved by a counting argument that has an inherent limitation, whereas Theorem \ref{rectLBthm} is proved by using flag algebras.
Theorems \ref{corAn1n2n3},  \ref{LBthm}, and \ref{rectLBthm}, provide  evidence for the next two conjectures.
 
\begin{conj} $\rcrn(K_{n_1,n_2,n_3})=A(n_1,n_2,n_3)$.
\end{conj}


\begin{conj} $\rcrn(K_{n_1,n_2,n_3})=\crn(K_{n_1,n_2,n_3})$.
\end{conj}

These two conjectures (if true) imply $\crn(K_{n_1,n_2,n_3})=A(n_1,n_2,n_3)$.\ms
  
A  complete multipartite graph is {\em balanced} if the partite sets all have the same cardinality.
In \cite{RT97MAA} it is shown that   $\lim_{n\to\infty}\frac {\crn(K_n)}{{n\choose 4}}\le \frac 3 8$ and $\lim_{n\to\infty}\frac{\crn(K_{n,n})}{{n\choose 2}^2}\le \frac 1 4$ and the limits exist.  We establish an analogous upper bound for the balanced complete $r$-partite graph.  The {\em maximum crossing number} of a graph $G$ is 
\[\CRN(G):=\max\{\crn(D) : D \mbox{ is a good drawing of }G\}.\]
With this notation,  it is shown in \cite{RT97MAA} that 
\[\lim_{n\to\infty}\frac {\crn(K_n)}{ \CRN(K_n)}\le \lim_{n\to\infty}\frac{H(n)}{\CRN(K_n)}=\frac 3 8\mbox{ and }\lim_{n\to\infty}\frac {\crn(K_{n,n})}{ \CRN(K_{n,n})}\le \lim_{n\to\infty}\frac{Z(n,n)}{\CRN(K_{n.n})}=\frac 1 4.\]

To state our bound for the complete multipartite graph, we need additional notation. The balanced complete $r$-partite graph $K_{n,\dots,n}$ will be denoted by
$\Krn r n$ because it is the join of $r$ copies of the complement of ${K_n}$.  Note that $\Krn 2 n=K_{n,n}$, $\Krn 3 n=K_{n,n,n}$, and $\Krn n 1=K_n$.

\begin{rem} {\rm The maximum crossing number can be computed as the number of choices of $4$ endpoints that can produce a crossing, and can be realized by a rectilinear drawing with vertices evenly spaced on a circle and vertices in the same partite set consecutive (this is well-known for the complete graph and complete bipartite graph).
Thus $ \CRN(K_n)={n \choose 4}$ and 
\beq \CRN(\Krn r n)={r \choose 2}{n \choose 2}^2+{r}{r-1 \choose 2}{n \choose 2}{n \choose 1}^2+{r \choose 4}{n \choose 1}^4,\label{CReq}\eeq
 with \eqref{CReq} obtained by choosing points partitioned among the partite sets as (2,2), (2,1,1), and (1,1,1,1).
For $r=2,3,4$ this yields
\ben
\item $ \CRN(\Krn 2 n)={n \choose 2}^2$,
\item $ \CRN(\Krn 3 n)=3{n \choose 2}^2+3{n \choose 2}{n \choose 1}^2$,
\item $ \CRN(\Krn 4 n)=6{n \choose 2}^2+12{n \choose 2}{n \choose 1}^2+{n \choose 1}^4$.
\een
 }
\end{rem}  

A {\em geodesic spherical drawing} of $G$ is a good drawing of $G$ obtained by placing the vertices of $G$  on  a sphere, drawing edges as geodesics, and projecting onto the plane. 
In a {\em random} geodesic drawing, the vertices are placed randomly on the sphere.  For integers $r\ge 2$ and $n\ge 1$,  define $s(r,n)$ to be the expected number of crossings in a random geodesic spherical drawing of $\Krn r n$ and  $\zeta(r):=\frac{3(r^2-r) }{8 \left(r^2+r-3\right)}$.  The next theorem is proved  in Section \ref{szeta}.

\begin{thm}\label{zetathm}  For $r\ge 2$, $\lim_{n\to\infty}\frac{s(r,n)}{\CRN(\Krn r n)}=\zeta(r)$.
\end{thm}

\begin{cor}  $\lim\sup_{n\to\infty}\frac{\crn(\Krn r n)}{\CRN(\Krn r n)}\le \zeta(r)$.
\end{cor}

\begin{obs} Note that $\zeta(r)=\frac{3(r^2-r) }{8 \left(r^2+r-3\right)}$ is monotonically increasing for $r\ge 3$, so  $\frac 1 4 =\zeta(2)=\zeta(3)<\zeta(4)<\cdots< \zeta(r)<\zeta(r+1)< \cdots< \frac 3 8$, and $\lim_{r\to\infty}\zeta(r)=\frac 3 8$.
\end{obs}

\begin{obs}  As $n\to\infty$, $\CRN(\Krn 3 n)\approx 3\frac {n^4} 4+ 3 \frac {n^4} 2 =\frac 9 4 n^4$ and $A(n,n, n)\approx 3 \left(\frac {n^4} {16}+\frac {n^4} {8}\right)=\frac 9 {16} n^4$, so $\lim_{n\to\infty} \frac{A(n,n,n)}{\CRN(\Krn 3 n)}=\frac 1 4=\zeta(3)$.
\end{obs}


\section{Proofs of Theorems \ref{corAn1n2n3}, \ref{LBthm}, and \ref{rectLBthm}}\label{sA3n}   In this section we  define a drawing of $K_{n_1,n_2,n_3}$ and use it to show  that $\rcrn(K_{n_1,n_2,n_3})\le A(n_1,n_2,n_2)$ for all $n_1,n_2,n_3$.  We also prove that asymptotically $0.666 A(n,n,n)\le \crn(K_{n,n,n})$ and $0.973 A(n,n,n) \le \rcrn(K_{n,n,n})$ for large $n$.

The standard way of producing a rectilinear drawing of the complete bipartite graph $K_{n,m}$ with $Z(n,m)$ crossings is a {\em $2$-line} drawing, constructed by drawing two perpendicular lines and placing the vertices of each partite set on one of the lines, with  about half of the points  on either side of the intersection of the lines. In the next definition we extend the idea of a 2-line drawing.    
\begin{figure}[h!]
\begin{center}
\includegraphics[scale=1.2]{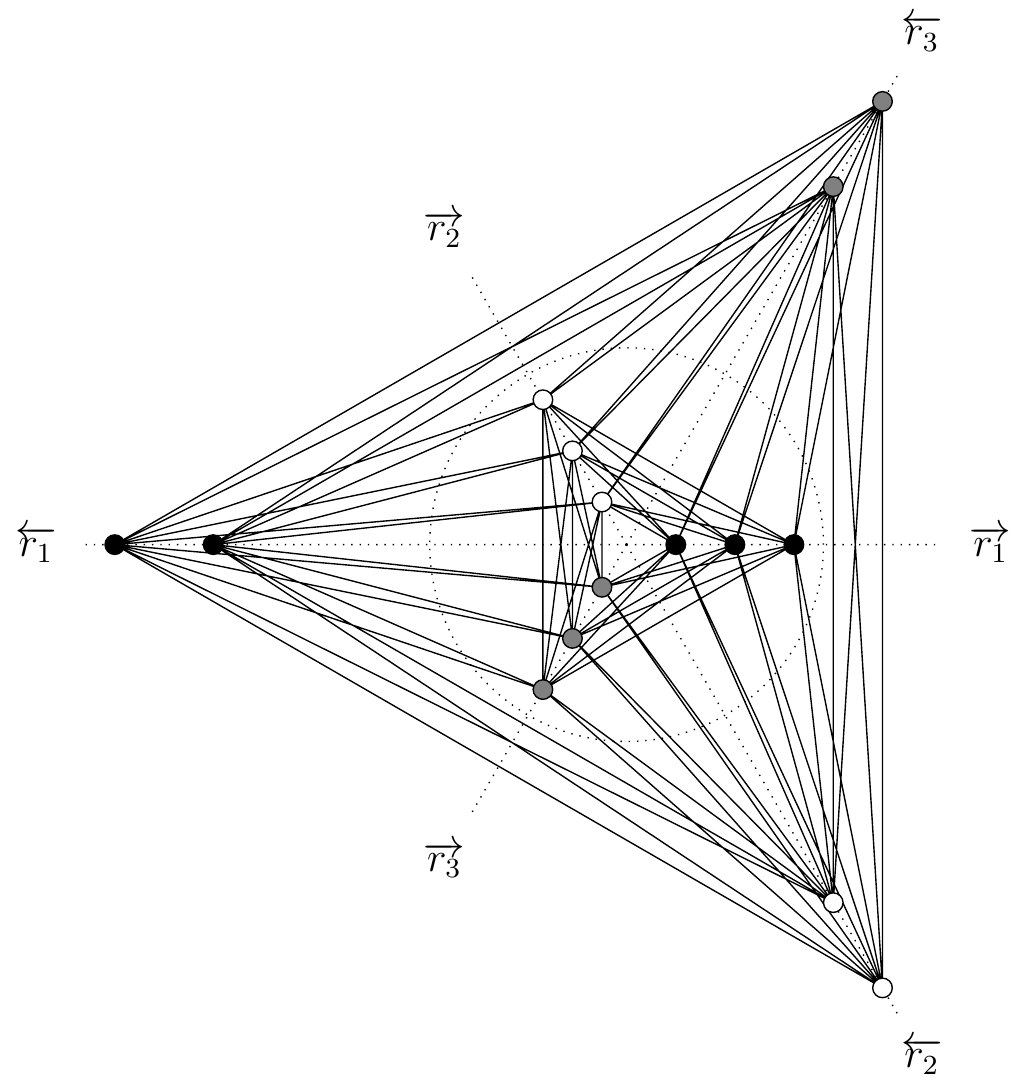}
\caption{An alternating 3-line drawing of $K_{5,5,5}$.  Points on opposite rays are in one partite set. The partite sets are distinguished by the color of the nodes.
The distances were slightly adjusted for visual clarity.  The rays and unit circle are shown faint and dotted.}
\label{alt3lineK555}\vspace{-15pt}
\end{center}
\end{figure}

\begin{defn}\label{alt3linedef} {\rm An {\em alternating $3$-line drawing} of $K_{n_1,n_2,n_3}$ is produced as follows:
\ben
\item\label{s1} Draw 3 rays $\overrightarrow{r_1},\overrightarrow{r_2},\overrightarrow{r_3}$  (called the {\em large rays}) 
that all originate from one point (called the {\em center}) with an angle of $120^\circ$ between each  pair of rays. 
\item For every $i \in \{1,2,3\}$, draw a ray $\overleftarrow{r_i}$ (called a {\em small ray})
from the center in the opposite direction of $\overrightarrow{r_1}$. We call $\overleftarrow{r_i}$  and $\overrightarrow{r_i}$ \emph{opposite} rays, and together they form the $i$th line $\ell_i$.
\item\label{spoints} For $i=1,2,3$:
\ben
\item Define  $a_i:=\left\lceil \frac{n_i}{2}\right\rceil$ and $b_i:=\left\lfloor \frac{n_i}{2}\right\rfloor $
\item\label{sa}  On $\overrightarrow{r_i}$, place $a_i$ points  at distances
$\frac 1 {a_i+1}, \frac 2 {a_i+1}, \dots \frac{a_i}{a_i+1}$ from the center.
\item\label{sb} On $\overleftarrow{r_i}$,  place  $b_i$ points at distances $3,4, \dots, b_i+2$ from the center. 
\een
\item\label{s6} For each pair of points not on the same line $\ell_i$, draw the line segment between the points.
\een
}\end{defn}

 The rays in Definition \ref{alt3linedef}  are not part of the drawing  but are useful reference terms.  Figure \ref{alt3lineK555} shows an alternating 3-line drawing of $K_{5,5,5}$. 


  The function  defined in \eqref{ALeq} below more naturally captures the number of crossings in an alternating 3-line drawing.
  
\begin{align}\label{ALeq} 
A_{3L}(n_1,n_2,n_3):=&\nonumber\\
\sum_{\substack{i=1,2,3\\ \{j,k\}=\{1,2,3\}\setminus \{i\}}}&\left[
{\left\lceil \frac{n_j}{2}\right\rceil \choose 2}{\left\lceil \frac{n_k}{2}\right\rceil \choose 2} +{\left\lceil \frac{n_j}{2}\right\rceil \choose 2}{\left\lfloor \frac{n_k}{2}\right\rfloor \choose 2} + {\left\lfloor \frac{n_j}{2}\right\rfloor \choose 2}{\left\lceil \frac{n_k}{2}\right\rceil \choose 2} + {\left\lfloor \frac{n_j}{2}\right\rfloor \choose 2}{\left\lfloor \frac{n_k}{2}\right\rfloor \choose 2} + \right.\nonumber\\
&\left.\left({\left\lceil \frac{n_i}{2}\right\rceil \choose 2}+{\left\lfloor \frac{n_i}{2}\right\rfloor \choose 2} \right)
\left(\left\lfloor \frac{n_j }{2}\right\rfloor\left\lceil \frac{n_k}{2}\right\rceil+\left\lceil \frac{n_j}{2}\right\rceil \left\lfloor \frac{n_k }{2}\right\rfloor \right)\right].
\end{align}

\begin{thm}\label{A3nThm} For $n_1,n_2,n_3 \ge 1$, an alternating $3$-line drawing of $K_{n_1,n_2,n_3}$ has at most\break$A_{3L}(n_1,n_2,n_3)$ crossings.
\end{thm} 
\bpf We count the maximum number of possible crossings in an alternating 3-line drawing of $K_{n_1,n_2,n_3}$.  There are two types of pairs of points that can result in crossings, (2,2) and (2,1,1), arising from choosing points partitioned among the partite sets as (2,2) and (2,1,1). Throughout this proof, $\{i,j,k\} = \{1,2,3\}$.

 For type (2,2),  if at least one pair of points has a point from each of the large and small ray,  i.e., from two opposite rays, we do not get a crossing. Thus we assume each set of two points in the same partite set is actually on the same ray.   We can choose any two rays that are not opposite, with each ray to contain two points.  There are 12 pairs of rays (omitting the opposite pairs), including  3 cases of $\left\lceil \frac{n_j}{2}\right\rceil$ and $\left\lceil \frac{n_k}{2}\right\rceil$ points,   3 cases of $\left\lfloor \frac{n_j}{2}\right\rfloor $ and $\left\lfloor \frac{n_k}{2}\right\rfloor$ points, and  6 cases of $\left\lceil \frac{n_j}{2}\right\rceil$ points and   $\left\lfloor \frac{n_k}{2}\right\rfloor $ points. 
Thus  there are at most 
\[ 
\sum_{\substack{i=1,2,3\\ \{j,k\}=\{1,2,3\}\setminus \{i\}}}\left[
{\left\lceil \frac{n_j}{2}\right\rceil \choose 2}{\left\lceil \frac{n_k}{2}\right\rceil \choose 2} +{\left\lceil \frac{n_j}{2}\right\rceil \choose 2}{\left\lfloor \frac{n_k}{2}\right\rfloor \choose 2} + {\left\lfloor \frac{n_j}{2}\right\rfloor \choose 2}{\left\lceil \frac{n_k}{2}\right\rceil \choose 2} + {\left\lfloor \frac{n_j}{2}\right\rfloor \choose 2}{\left\lfloor \frac{n_k}{2}\right\rfloor \choose 2}\right]
\] crossings of type (2,2).

 Consider pairs of points partitioned as type (2,1,1).  
 Denote by $B$ the unit ball centered at the center of the drawing.
 Observe that the line segment between any two points on small rays is disjoint from $B$ and the line segment between any two points on large rays is  entirely in $B$.
 If we choose the two points in the same partite set from opposite rays, then we do not get a crossing. Thus we assume the   two points in the same partite set are actually on the same ray.

We can choose any one ray to contain the two points from a (2,1,1) partition of points.  
Suppose the ray chosen is $\overrightarrow{r_i}$ or $\overleftarrow{r_i}$, where $i \in \{1,2,3\}$.
Line $\ell_i$ containing $\overrightarrow{r_i}$ and $\overleftarrow{r_i}$ divides the plane to two half-planes. 
To have a crossing, the other two points must come from the same half-plane.
Thus the number of choices of pairs of points  from the two rays in one half plane is 
$\left\lfloor \frac{n_j }{2}\right\rfloor\left\lceil \frac{n_k}{2}\right\rceil+\left\lceil \frac{n_j}{2}\right\rceil \left\lfloor \frac{n_k }{2}\right\rfloor$. 
Each of these is multiplied by the choice of pair from $\overrightarrow{r_i}$ and $\overleftarrow{r_i}$,
which gives the maximum number of crossings containing a pair from $\overrightarrow{r_i}$ or $\overleftarrow{r_i}$ as 
\beq
\left({\left\lceil \frac{n_i}{2}\right\rceil \choose 2}+{\left\lfloor \frac{n_i}{2}\right\rfloor \choose 2} \right)
\left(\left\lfloor \frac{n_j }{2}\right\rfloor\left\lceil \frac{n_k}{2}\right\rceil+\left\lceil \frac{n_j}{2}\right\rceil \left\lfloor \frac{n_k }{2}\right\rfloor \right).\label{one211}\vspace{-4pt}
\eeq
The maximum number of crossings of type $(2,1,1)$ is obtained by summing \eqref{one211} over all choices of $i \in \{1,2,3\}$. 

Thus the total number of crossings in this drawing is at most $A_{3L}(n_1,n_2,n_3)$. 
 \epf

To complete the proof of Theorem \ref{corAn1n2n3}, we show that $A_{3L}(n_1,n_2,n_3)= A(n_1,n_2,n_3)$. 
First we show that for all $a$  we have ${ \lceil \frac{a}{2} \rceil \choose 2} + { \lfloor \frac{a}{2} \rfloor \choose 2} = \lfloor \frac{a}{2} \rfloor \lfloor \frac{a-1}{2} \rfloor$.
By distinguishing odd and even case we get
\[
{ \lceil \frac{a}{2} \rceil \choose 2} + { \lfloor \frac{a}{2} \rfloor \choose 2} = \begin{cases}  
  2{\frac{a}{2} \choose 2} = \frac{a}{2}(\frac{a}{2} -1) =  \lfloor \frac{a}{2} \rfloor \lfloor \frac{a-1}{2} \rfloor &   \text{ for } a \text { even,} \\
 \frac{1}{2} \left( \frac{a+1}{2} \right) \left( \frac{a+1}{2} - 1 \right) + \frac{1}{2} \left( \frac{a-1}{2} \right) \left( \frac{a-1}{2} - 1 \right) =  \frac{a^2 - 2a + 1}{4}  =  \lfloor \frac{a}{2} \rfloor \lfloor \frac{a-1}{2} \rfloor   &   \text{ for } a \text { odd.}
\end{cases}
\]

Next we show $\lfloor \frac{a}{2} \rfloor \lceil \frac{b}{2} \rceil + \lceil \frac{a}{2} \rceil \lfloor \frac{b}{2} \rfloor  = \lfloor \frac{ab}{2} \rfloor$.
Again, we distinguish cases by the parity of $a$ and $b$ and obtain
\[
\bigg\lfloor \frac{a}{2} \bigg\rfloor \bigg\lceil \frac{b}{2} \bigg\rceil + \bigg\lceil \frac{a}{2} \bigg\rceil \bigg\lfloor \frac{b}{2} \bigg\rfloor =
\begin{cases}
\frac{a}{2} \left( \lceil \frac{b}{2} \rceil + \lfloor \frac{b}{2} \rfloor \right) = \frac{ab}{2} = \lfloor \frac{ab}{2} \rfloor & \text{ for } a \text{ even},\\
\frac{(a-1)(b+1)}{4} + \frac{(a+1)(b-1)}{4} = \frac{2ab -2}{4} = \lfloor \frac{ab}{2} \rfloor  & \text{ for } a,b \text{ odd.}\\
\end{cases}
\]

Using these two observations it is straightforward to show $A(n_1,n_2,n_3)=A_{3L}(n_1,n_2,n_3)$.
The assertion that $\rcrn(K_{n_1,n_2,n_3})\le A(n_1,n_2,n_3)$ then follows from Theorem \ref{A3nThm}.  This completes the proof of Theorem \ref{corAn1n2n3}.\ms

Next we prove Theorem \ref{LBthm}, i.e., $0.666 A(n,n,n) \le \crn(K_{n,n,n})$ for large $n$.

\bpf It is known  that $\crn(K_{2,3,n})=4\left\lfloor \frac{n}{2}\right\rfloor \left\lfloor \frac{n-1}{2}\right\rfloor + n$ (see \cite{A86}).  So each copy of $K_{2,3,n}$ in $K_{n,n,n}$ has approximately $ n^2$ crossings.  The number of copies of $K_{2,3,n}$ in $K_{n,n,n}$ is $6{n\choose 2}{n\choose 3} {n\choose n}$, where the factor of 6 comes from choosing which of the  three partite sets in $K_{n,n,n}$ is used for the 2, which for the 3, and which for the $n$.  Thus we count about
$ 
\left(  n^2\right)\left(6 \cdot \frac {n^2} 2\cdot\frac {n^3} 6 \right)=\frac 1 {2} n^7
$ crossings (counting each crossing multiple times).

The number of times a crossing gets counted varies with whether the end points are partitioned of type (2,2) or type (2,1,1).  For type (2,2), we can arrange the $K_{2,3,n}$ among the three partite sets   as (2,2,0), (2,0,2), or type (0,2,2), and in each case there are 2 choices.  Thus a crossing of type (2,2) is counted 
\bea
2\left[{{n-2}\choose{0}}{{n-2}\choose{1}}{{n}\choose{n}} +
{{n-2}\choose{0}}{{n}\choose{3}}{{n-2}\choose{n-2}} + 
{{n}\choose{2}}{{n-2}\choose{1}}{{n-2}\choose{n-2}}\right]& \approx &\\
2\left[n + \frac{n^3}{6}+ \frac{n^3}{2}\right]&\approx &  \frac{4n^4}{3}\vspace{-5pt}
\eea
times.

For type (2,1,1), we can arrange the $K_{2,3,n}$ among the three partite sets   as (2,1,1), (1,2,1), or type (1,1,2), and in each case there are 2 choices.  Thus a crossing of type (2,1,1) is counted 
\bea
2\left[{{n-2}\choose{0}}{{n-1}\choose{2}}{{n-1}\choose{n-1}} \!+\!
{{n-1}\choose{1}}{{n-2}\choose{1}}{{n-1}\choose{n-1}} \!+\! 
{{n-1}\choose{1}}{{n-1}\choose{2}}{{n-2}\choose{n-2}}\right]\!\!\!\!& \approx &\!\!\!\!\null\\
2\left[\frac {n^2} 2 + n^2+ \frac{n^3}{2}\right]&\!\!\approx & \!\!\!\! n^3\vspace{-5pt}
\eea
times.

Since $\frac{4n^3}{3}>n^3$ and $A(n,n,n)\approx \frac 9 {16} n^4$, asymptotically we have at least 
\[\frac{\frac 1 2 n^7}{\frac{4n^3}{3}}=\frac 3 8 n^4=\frac 2 3 \left(\frac 9 {16} n^4\right) = \frac 2 3 A(n,n,n)> 0.666 A(n,n,n)\]
crossings.
\epf

\begin{rem}{\rm We point out that the counting method  used in the proof of Theorem \ref{LBthm} has a structural limitation.  We use the count number for a (2,2) partition  as the number of times a $K_{2,3,n}$ is counted (because it is the larger), even though we know that asymptotically 2/3 of the crossings in an alternating 3-line drawing of $K_{n,n,n}$ are of type (2,1,1) rather than (2,2).  So even with the assumption that $\crn(K_{n,n,n})=A(n,n,n)$,  this method cannot be expected to produce a lower bound of $cA(n,n,n)$ with $c$ close to 1.
}\end{rem}

Finally we prove Theorem \ref{rectLBthm}, i.e., $0.973 A(n,n,n) \le \rcrn(K_{n,n,n})$ for large $n$.  The proof uses flag algebras, a method developed by Razborov~\cite{Raz07}.  A brief explanation of this technique specific to its use in our proof can also be found in Appendix \ref{appdx1}.
We use an approach similar to the technique Norin and Zwols~\cite{NorinZwols} used to show that $
0.905 Z(m,n) \leq \crn(K_{m,n})$; however, we restrict our attention  to rectilinear drawings.

\bpf For sufficiently large $n$, we first use flag algebra to methods show that in any  rectilinear drawing of $K_{n,n,n}$ the average number of crossings over all the copies of $K_{3,2,2}$ that appear in $K_{n,n,n}$ is greater than 5.6767.  
In our application, we record in flags crossings and tripartitions. We ignore rest of the embedding.

Let $G$ be a tripartite graph on $n$ vertices with a rectilinear drawing.
A corresponding flag $F_G$ on $n$ vertices $V$ contains a function $\varrho_1: V \rightarrow \{0,1,2\}$,
which records the partition of the vertices, and a function $\varrho_2: V^4 \rightarrow \{0,1\}$,
which record crossings. We define $\varrho_2(a_1,a_2,b_1,b_2) = 1$ if the vertices of $G$
corresponding to $a_1$ and $a_2$ form an edge of $G$ that crosses an edge of $G$ formed by $b_1$ and $b_2$, and 0 otherwise.

We use flags on 7 vertices  obtained from rectilinear drawings of $K_{3,2,2}$ (so $m=7$   in Equation \eqref{eq:sum} in Appendix  \ref{appdx1}).
All rectilinear drawings of $K_7$ were obtained by Aichholzer, Aurenhammer and Krasser~\cite{AichholzerAK:2002}.
The drawings give us 6595 flags. 
We generate 42 types, which leads to 42 equations like \eqref{eq:theta} in Appendix  \ref{appdx1}.
The optimal linear combination of these equations is computed by CSDP~\cite{Borchers:1999}, an open source semidefinite
program solver.
CSDP is a numerical solver that provides a positive semidefinite matrix $M$ of floating point numbers.
We round the matrix $M$ in Sage~\cite{sw:sage} to a positive semidefinite matrix $Q$ with rational entries.  The rounding is done by decomposing $M = U^T D U$ (where $D$ is a diagonal matrix of eigenvalues and $U$ is a real orthogonal matrix of eigenvectors),  rounding the entries of $D$ and $U$ to rational matrices $\hat D$ and $\hat U$, and constructing matrix $Q=\hat U^T \hat D\hat U$.  
Then we use $Q$ to compute the resulting bound 
$1419186177261/250000000000 > 5.6767$.
Software needed to perform the whole computation is available at \oururl.

In a complete graph on 7 vertices, the number of 4-tuples of points is ${7 \choose 4}=35$.  Thus the `density' of crossings in $K_{3,2,2}$ is at least $\frac {5.6767} {35}$.  The graph $K_{n,n,n}$ must have at least this density times the number of 4-tuples, and the number of 4-tuples is ${3n\choose 4}\approx \frac {81n^4}{24}$.  Since $A(n,n,n)\approx \frac {9n^4}{16}$, asymptotically
\[\rcrn(K_{n,n,n})> \left(\frac {5.6767} {35}\right)\left(\frac {81n^4}{24}\right)
\approx \left(\frac {5.6767} {35}\right)6A(n,n,n)>.973A(n,n,n).\qedhere\]
\epf

\begin{rem} {\rm The flag algebra method just applied to $\rcrn(K_{n,n,n})$ with $n\to\infty$ 
will also work for $\rcrn(K_{n_1,n_2,n_3})$ where $n_i\to\infty$ for all $i=1,2,3$.}\end{rem} 

\section{Proof of Theorem \ref{zetathm}}\label{szeta}

We need a preliminary lemma.

  \begin{lem}\label{paircrossprob}  In a random geodesic spherical drawing of a pair of disjoint edges, the probability that the pair  crosses is $\frac 1 8$.
\end{lem}
\bpf  A pair of edges is determined by two sets of endpoints.  Each set of two endpoints determines a great circle, and these two great circles  intersect in two antipodal points.  These two antipodal points of intersection are the potential crossing points, and a crossing occurs if and only if both edges include the same antipodal point.  Notice that first picking two great circles uniformly at random, and then picking two points uniformly at random from each of the great circles is equivalent to picking two pairs of points uniformly at random from the sphere. Therefore, for each set of two endpoints, the probability that the great circle geodesic between them includes one of the two antipodal points is $\frac 1 2$, so the probability that both edges include an antipodal point is $\frac 1 4$.  Half the time these are the same antipodal point.    \epf

We are now ready to prove Theorem \ref{zetathm}, i.e., $\lim_{n\to\infty}\frac{s(r,n)}{\CRN(\Krn r n)}=\zeta(r)$.

\bpf
The probability of getting a crossing among four points in a geodesic spherical drawing of $\Krn r n$ depends on how the points are partitioned among the partite sets, because different partitions of four points  have different numbers of  pairs of disjoint edges.  Define three types of partitions of four points, classified by the number of pairs (of disjoint edges) produced. 

\noi{\bf Type A} 0 pairs:  The four points are  partitioned among partite sets as (4) or (3,1).  
Let $\alpha_r$ denote the probability that four randomly chosen points in $\Krn r n$ are of this type.

\noi{\bf Type B} 2 pairs: The four points are  partitioned among partite sets as (2,2) or (2,1,1). Let $\beta_r$ denote the probability that four randomly chosen points in $\Krn r n$ are of this type.

\noi{\bf Type C} 3 pairs: The four points are  partitioned among partite sets as (1,1,1,1). Let $\gamma_r$ denote the probability that four randomly chosen points in $\Krn r n$ are of this type.

We assume that $n$ is large relative to $r$, so we can ignore the difference between $n-1$ and $n$, etc., and we  focus only on which partite sets are chosen.
For Type C we must choose four distinct partite sets, so $\gamma_r=\frac {r(r-1)(r-2)(r-3)}{r^4}=\frac {(r-1)(r-2)(r-3)}{r^3}$.  For Type A there are two choices.  For partition (4) the probability is $\frac 1 {r^3}$.  To determine the probability of partition (3,1) we count the ways that we can choose 4 partite sets with 3 of them being the same set (which we call a (3,1) choice), and divide by $r^4=$ the number of all possible arrangements of four points into $r$ partite sets. A  (3,1) choice can be made  by first choosing two distinct partite sets (there are $r(r-1)$ ways to select the two, with the first choice to appear 3 times) and then indicating the order of these partite sets  (there are four different orders, determined by where  the singleton is placed in the order).  So the probability is $\frac {4r(r-1)} {r^4}=\frac {4(r-1)} {r^3}$.  Thus $\alpha_r=\frac {4(r-1)} {r^3}+\frac {1} {r^3}=\frac {4r-3} {r^3}$.  Then $\beta_r=1-\alpha_r-\gamma_r$.  

Let $q$ be the number of 4-tuples of points.  By Lemma \ref{paircrossprob}, the expected number of crossings in a geodesic spherical drawing is $\frac 1 8$ the number of  pairs of disjoint edges, and the number of pairs  is 
\[(3\gamma_r+2\beta_r)q=(3\gamma_r+2(1-\alpha_r-\gamma_r))q=(2+\gamma_r-2\alpha_r)q,\]  
so $s(r,n)=\frac 1 8 (2+\gamma_r-2\alpha_r)q$.
 In the earlier described drawing that maximizes the number of crossings, every 4-tuple of Type B and C produces one crossing.  There are $(\beta_r+\gamma_r)q$ such 4-tuples, and therefore
\[\CRN\left(\Krn r n\right)=(\gamma_r+\beta_r)q
=(1-\alpha_r)q.\]
  Thus
\bea \lim_{n\to\infty}\frac{s(r,n)}{\CRN(\Krn r n)}&=& \frac{\frac 1 8 (2+\gamma_r-2\alpha_r)q}{(1-\alpha_r)q}\\
&=&\frac{2 r^3 + (r - 1) (r - 2) (r - 3) - 2 (4 r - 3)}{8
   \left(r^3 - (4 r - 3)\right)}\\
   &=&\frac{3(r^2-r) }{8 \left(r^2+r-3\right)}\\
   &=&\zeta(r).\hspace{6cm}\qedhere\eea
\epf

\subsection*{Acknowledgments}
This research began at the American Institute of Mathematics workshop Exact Crossing Numbers, and the authors thank AIM.  The authors thank Sergey Norin for many helpful conversations during that workshop.  This paper was finished while Hogben, Lidick\'y, and Young were general members in residence at the Institute for Mathematics and its Applications, and they thank IMA.   The authors also thank NSF for their support of these institutes. The work of Gethner and Pfender is supported in part by their respective
Simons Foundation Collaboration Grants for Mathematicians.  The work of Lidick\'y is partially supported by NSF grant DMS-12660166.

\appendix

\section{Flag Algebras}\label{appdx1}

The theory of flag algebras is a recent framework developed by Razborov~\cite{Raz07}.
The method was designed to attack Tur\'an and subgraph density problems in extremal combinatorics and has been applied to graphs~\cite{HatamiJKNR:2012}, hypergraphs~\cite{Falgas2013},
geometry~\cite{KralMachSereni2012}, permutations~\cite{permutations}, and crossing numbers~\cite{NorinZwols}, to name some.
For more applications see a recent survey by Razborov~\cite{Razborov13flagalgebras}.

Use of flag algebra methods usually depends on a computer program that
generates a large semidefinite program, which can be solved
by an available solver. The method is in some cases automated by
Flagmatic \cite{flagmatic}. However, Flagmatic does not support counting
crossings. Hence we developed our own software, available at \oururl; our use of this software  is described in the proof of Theorem  \ref{rectLBthm}.

Rather than attempt to give a formal setup of the framework of flag algebras, this introduction is intended to give the reader enough background  to understand how we apply the method to prove Theorem 1.3. 
For a formal description of the method, involving the algebra of linear combinations of non-negative homomorphisms, see Rasborov~\cite{Razborov13intro}.

\subsection{Densities}\label{ssdense}

Let $G$ be a large graph on $n$ vertices and let $d_P(G)$ be the density of property  $P$ in $G$. In our case, the property $P$ is a crossing. We can compute $d_P(G)$ by computing $d_P(H)$ for all possible graphs $H$ on $m$ vertices, where $m << n$, and then count how often $H$ appears in $G$. We denote the density of $H$ in $G$
by $d_H(G)$, which is the same as the probability that $m$ vertices of $G$ selected uniformly at random induce a copy of $H$.
This gives the following equality, 
\beq
d_P(G) = \sum_{|V(H)|=m}  d_P(H)d_H(G).\label{eq:sum}
\eeq
Therefore, depending on how we are optimizing, we attain one of the following inequalities, 
\[
\min_{|V(H)|=m} d_P(H) \le d_P(G) \le \max_{|V(H)|=m} d_P(H). 
\]

In general, these bounds tend to be rather weak, so flag algebras are used to improve inequalities on $d_H(G)$. Assuming there exists the linear inequality 
\[0 \le \displaystyle\sum_{|V(H)|=m} c_H d_H(G),\] then 
\[d_P(G) \ge \sum_{|V(H)|=m} (d_P(H) - c_H) d_H(G) \ge \min_{|V(H)|=m} (d_P(H) - c_H).\] This may improve the bound if there are negative values for $c_H$. Semidefinite programming is used to determine these coefficients.

\subsection{Flags}\label{ssflag}

A \emph{type} $\sigma$ is a graph on $s$ vertices with a bijective labeling function $\theta: [s] \to V(\sigma)$. A $\sigma$-\emph{flag} $F$ is a graph $H$ containing an induced copy of $\sigma$ labeled by $\theta$ and  the order of $F$ is $ |V(F)|$. 
Let $\ell$, $m$, and $s$ be integers such that $s < m$, and $2\ell \le m + s$. These values ensure that a graph on $m$ vertices can have $\sigma$-flags of order $\ell$ that intersect in exactly $s$ vertices. Define $\mathcal{F}_\ell^{\sigma}$ to be the set of all $\sigma$-flags on $\ell$ vertices, up to isomorphism.

Given an injection from $[s] \to V(G)$, $\theta$ and $F \in \mathcal{F}_\ell^{\sigma}$, define $d_{F}(G; \theta)$ to be the density of $F$ in $G$ labeled by $\theta$. Note that for $|\sigma| = 0$, this density corresponds to the $d_F(G)$. If $F_a, F_b \in \mathcal{F}_\ell^{\sigma}$, then we say $d_{F_a,F_b}(G; \theta)$ is the density of the graph created when $F_a$ and $F_b$ intersect exactly at $\sigma$.

\begin{thm}[Razborov~\cite{Raz07}]\label{raz} 
For any $F_a, F_b \in \mathcal{F}_\ell^{\sigma}$ and $\theta$, \[d_{F_a}(G; \theta) d_{F_b}(G; \theta) = d_{F_a,F_b}(G; \theta) + o(1).\]
\end{thm}


Let $\ff$ be a vector with entries $d_{F_i}(G; \theta)$ for all $F_i \in \mathcal{F}_\ell^{\sigma}$ and let $Q$ be a positive semidefinite matrix with $q_{ij}$ as the $ij$th entry. Then we get \[0 \le \ff^T Q \ff  =\sum_{F_i,F_j \in \mathcal{F}_\ell^{\sigma}} q_{ij}  d_{F_i}(G; \theta)  d_{F_j}(G; \theta).\] Theorem \ref{raz} gives  \[0 \le  \displaystyle\sum_{F_i,F_j \in \mathcal{F}_\ell^{\sigma}} q_{ij}  d_{F_i,F_j}(G; \theta) + o(1).\]

By averaging over all $\theta$ and all subgraphs on $m$ vertices, we can obtain an inequality of the form  
\beq
0 \le \displaystyle\sum_{H \in \mathcal{F}_m^{0}} c_H d_H(G)+o(1),\label{eq:theta}
\eeq
where $c_H$ is a function of $\sigma$, $m$, and $Q$. So asymptotically as $n\to \infty$,
\[d_P(G) \ge \sum_{H \in \mathcal{F}_m^{0}} (d_P(H) - c_H) d_H(G) + o(1) \ge \min_{H \in \mathcal{F}_\ell^{0}} (d_P(H) - c_H).\]



\begin{thebibliography}{99}

\bibitem{AichholzerAK:2002}
O.~Aichholzer, F.~Aurenhammer, and H.~Krasser.
\newblock Enumerating order types for small point sets with applications.
\newblock {\em Order}, 19: 265--281, 2002.

\bibitem{A86}
 K. Asano. The crossing number of $K_{1,3,n}$ and $K_{2,3,n}$. 
{\em  Journal of Graph
Theory}, 10: 1--8, 1986. 

\bibitem{BW10} L. Beineke and R. Wilson.  The early history of the brick factory problem. {\em Mathematical  Intelligencer},  32: 41--48, 2010.





\bibitem{permutations}
J.~Balogh, P.~Hu, B.~Lidick\'y, O.~Pikhurko, B.~Udvari, and J.~Volec.
\newblock Minimum number of monotone subsequences of length 4 in permutations.
To appear in {\em Combinatorics, Probability and Computing}.

\bibitem{Borchers:1999}
B.~Borchers.
\newblock {CSDP, A C} library for semidefinite programming.
\newblock {\em Optimization Methods and Software}, 11: 613--623, 1999.

\bibitem{Falgas2013}
V.~Falgas-Ravry and E.~R. Vaughan.
\newblock Applications of the semi-definite method to the {T}ur\'an density
  problem for 3-graphs.
\newblock {\em Combinatorics, Probability and Computing}, 22: 21--54, 2013.

\bibitem{HatamiJKNR:2012}
H.~Hatami, J.~Hladk{\'y}, D.~Kr{\'a}\soft{l}, S.~Norine, and A.~Razborov.
\newblock Non-three-colourable common graphs exist.
\newblock {\em Combinatorics, Probability and Computing}, 21: 734--742, 2012.

\bibitem{Ho08}  P. T. Ho.  The crossing number of $K_{1,m,n}$.  {\em  Discrete Mathematics},  308: 5996--6002, 2008. 


\bibitem{Ho13}
 P. T. Ho.  The crossing number of $K_{2,4,n}$. 
{\em  Ars Combinatorica}, 109: 527--537, 2013. 

\bibitem{HZ08}
 Y. Huang and T. Zhao.  The crossing number of $K_{1,4,n}$. 
{\em  Discrete Mathematics},  308: 1634--1638, 2008. 

\bibitem{KralMachSereni2012}
D.~Kr{\'a}{\soft{l}}, L.~Mach, and J.-S. Sereni.
\newblock A new lower bound based on {G}romov's method of selecting heavily
  covered points.
\newblock {\em Discrete Computational Geometry}, 48: 487--498, 2012.


\bibitem{NorinZwols}
S. Norin.
Tur\'{a}n's brickyard problem and flag algebras.
 Talk at {\em Geometric and topological graph theory} workshop, Banff International Research Station,  2013.  Video available at \url{http://www.birs.ca/events/2013/5-day-workshops/13w5091/videos/watch/201310011538-Norin.mp4}.

\bibitem{Raz07}
A.~A. Razborov.
 Flag algebras.
 {\em Journal of Symbolic Logic}, 72: 1239--128, 2007.

\bibitem{Razborov13flagalgebras}
A.~A. Razborov.
 Flag algebras: an interim report, 2013. To appear in the {\em Erd\"os Centennial Volume}, preliminary
text available at \url{http://people.cs.uchicago.edu/~razborov/files/flag_survey.pdf}.

\bibitem{Razborov13intro}
A.~A. Razborov.
\newblock What is {$\ldots$} a flag algebra?
\newblock {\em Notices American Mathematical Society}, 60 :1324--1327, 2013.

\bibitem{RT97MAA}
 R. B. Richter and C. Thomassen.  Relations Between Crossing Numbers of Complete and Complete Bipartite Graphs.
{\em  The American Mathematical Monthly},  104: 131--137, 1997. 

\bibitem{sw:sage}
W.~Stein et~al.
\newblock {\em {S}age {M}athematics {S}oftware} ({V}ersion 6.4), 2014.
\newblock The Sage Development Team, \url{http://www.sagemath.org}.

\bibitem{Turan} P. Tur\'an. A note of welcome. {\em  Journal of Graph
Theory},  1: 7--9, 1977.

\bibitem{flagmatic} E. R. Vaughan. Flagmatic
software.  Available at \url{http://flagmatic.org/}.

\bibitem{Zarank} K. Zarankiewicz. On a problem of P. Turan concerning graphs. {\em Fundamenta Mathematicae} 41: 137--14, 1954.
\end{thebibliography}
\end{document}